\documentclass{article}
\usepackage{amsmath,amsfonts,amsthm}
\usepackage[all]{xy}

\newtheorem{defn}{Definition}[section]
\newtheorem{example}[defn]{Example}
\newtheorem{prop}[defn]{Proposition}
\newtheorem{cor}[defn]{Corollary}
\newtheorem{rem}[defn]{Remark}
\newtheorem{thm}[defn]{Theorem}
\newtheorem{lemma}[defn]{Lemma}

\newcommand{\rank}{\mathrm{rk}}
\newcommand{\Cohstack}{\mathcal{C}oh}
\newcommand{\Bunstack}{\mathcal{B}un}
\newcommand{\Stabstack}{\mathcal{S}tab}
\newcommand{\Semistabstack}{\mathcal{S}stab}

\newcommand{\Parstack}{\mathcal{P}ar}
\newcommand{\Extstack}{\mathcal{E}xt}

\newcommand{\Gr}{\mathrm{Gr}}
\newcommand{\GL}{\mathrm{GL}}
\newcommand{\PGL}{\mathrm{PGL}}
\newcommand{\cohom}{\mathrm{H}}
\newcommand{\Hom}{\mathrm{Hom}}
\newcommand{\RHom}{\mathrm{RHom}}
\newcommand{\Hombdl}{\mathit{Hom}}
\newcommand{\Ext}{\mathrm{Ext}}

\newcommand{\Endbdl}{End}
\newcommand{\Aut}{\mathrm{Aut}}
\newcommand{\longto}[1][]{\stackrel{#1}{\longrightarrow}}
\newcommand{\otgnol}[1][]{\stackrel{#1}{\longleftarrow}}
\newcommand{\Quot}{\mathrm{Quot}}
\newcommand{\Flag}{\mathrm{Flag}}

\newcommand{\Pic}{\mathrm{Pic}}
\newcommand{\stackM}{\mathcal{M}}
\newcommand{\coarseM}{\mathfrak{Bun}}
\newcommand{\coarseU}{\mathfrak{U}}
\newcommand{\coarseGr}{\mathfrak{Gr}}
\renewcommand{\O}{\mathcal{O}}
\newcommand{\univ}{\mathrm{univ}}
\newcommand{\stab}{\mathrm{stab}}
\newcommand{\Suniv}{\mathcal{S}^\univ}
\newcommand{\Euniv}{\mathcal{E}^\univ}
\newcommand{\Funiv}{\mathcal{F}^\univ}
\newcommand{\Eprimeuniv}{\mathcal{E}'^\univ}
\newcommand{\stackU}{\mathcal{U}}
\newcommand{\stackV}{\mathcal{V}}
\newcommand{\stackW}{\mathcal{W}}

\newcommand{\integers}{\mathbb{Z}}
\newcommand{\dual}{\mathrm{dual}}
\newcommand{\coker}{\mathrm{coker}}
\newcommand{\id}{\mathrm{id}}
\newcommand{\im}{\mathrm{im}}
\newcommand{\hcf}{\mathrm{hcf}}
\newcommand{\Gm}{\mathbb{G}_m}
\newcommand{\bbAut}{\mathbb{A}\mathrm{ut}}

\begin{document}
\bibliographystyle{plain}

\title{Moduli stacks of vector bundles on curves and the King--Schofield rationality proof}
\author{Norbert Hoffmann\thanks{Mathematisches Institut der Universit\"at G\"ottingen, Bunsenstr. 3--5, D-37073 G\"ottingen, Germany. email: hoffmann@uni-math.gwdg.de}
}
\date{}
\maketitle

\section*{Introduction}
Let $C$ be a connected smooth projective curve of genus $g \geq 2$ over an algebraically closed field $k$. Consider the coarse moduli scheme $\coarseM_{r, d}$
(resp. $\coarseM_{r, L}$) of stable vector bundles on $C$ with rank $r$ and degree $d \in \integers$ (resp. determinant isomorphic to the line bundle $L$ on $C$).

Motivated by work of A. Tyurin \cite{tyurin, tyurin2} and P. Newstead \cite{newstead, newstead2}, it has been believed for a long time that $\coarseM_{r, L}$ is rational
if $r$ and the degree of $L$ are coprime. Finally, this conjecture was proved in 1999 by A. King and A. Schofield \cite{king-schofield}; they deduce it from their
following main result:
\begin{thm}[King--Schofield] \label{introthm}
  $\coarseM_{r, d}$ is birational to the product of an affine space $\mathbb{A}^n$ and $\coarseM_{h, 0}$ where $h$ be the highest common factor of $r$ and $d$.
\end{thm}
The present text contains the complete proof of King and Schofield translated into the language of algebraic stacks. Following their strategy, the moduli stack
$\Bunstack_{r, d}$ of rank $r$, degree $d$ vector bundles is shown to be birational to a Grassmannian bundle over $\Bunstack_{r_1, d_1}$ for some $r_1 < r$; then induction
is used. However, this Grassmannian bundle is in some sense twisted. Mainly for that reason, King and Schofield need a stronger induction hypothesis than \ref{introthm}:
They add the condition that their birational map preserves a certain Brauer class $\psi_{r, d}$ on $\coarseM_{r, d}$. One main advantage of the stack language here is that
this extra condition is not needed: The stack analogue of theorem \ref{introthm} is proved by a direct induction.

(In more abstract terms, this can be understood roughly as follows: A Brauer class corresponds to a gerbe with band $\Gm$. But the gerbe on $\coarseM_{r, d}$ corresponding
to $\psi_{r, d}$ is just the moduli stack $\Bunstack_{r, d}$. Thus a rational map of coarse moduli schemes preserving this Brauer class corresponds to a rational map of the
moduli stacks.)

This paper consists of four parts. Section \ref{statement} contains the precise formulation of the stack analogue \ref{mainthm} to theorem \ref{introthm}; then the original
results of King and Schofield are deduced. Section \ref{grassmannians} deals with Grassmannian bundles over stacks because they are the main tool for the proof of theorem
\ref{mainthm} in section \ref{proof}. Finally, appendix \ref{stacks} summarizes the basic properties of the moduli stack $\Bunstack_{r, d}$ that we use. In particular,
a proof of Hirschowitz' theorem about the tensor product of general vector bundles on $C$ is given here, following Russo and Teixidor \cite{russo-teixidor}.

The present article has grown out of a talk in the joint seminar of U. Stuhler and Y. Tschinkel in G\"ottingen. I would like to thank them for the opportunity to speak and
for encouraging me to write this text. I would also like to thank J. Heinloth for some valuable suggestions and for many useful discussions about these stacks.

\section{The King-Schofield theorem in stack form} \label{statement}

We denote by $\Bunstack_{r, d}$ the moduli stack of vector bundles of rank $r$ and degree $d$ on our smooth projective curve $C$ of genus $g \geq 2$ over $k = \bar{k}$.
This stack is algebraic in the sense of Artin, smooth of dimension $(g-1)r^2$ over $k$ and irreducible; these properties are discussed in more detail in the appendix.

Our main subject here is the birational type of $\Bunstack_{r, d}$. We will frequently use the notion of a rational map between algebraic stacks; it is defined in the
usual way as an equivalence class of morphisms defined on dense open substacks. A birational map is a rational map that admits a two-sided inverse.
\begin{defn}
  A rational map of algebraic stacks $\xymatrix{\stackM \ar@{-->}[r] & \stackM'}$ is \emph{birationally linear} if it admits a factorization
  \begin{displaymath} \xymatrix{
    \stackM \ar@{-->}[r]^-{\sim} & \stackM' \times \mathbb{A}^n \ar[r]^-{\mathrm{pr}_1} & \stackM'
  } \end{displaymath}
  into a birational map followed by the projection onto the first factor.
\end{defn}

Now we can formulate the stack analogue of the King-Schofield theorem \ref{introthm}; its proof will be given in section \ref{proof}.
\begin{thm} \label{mainthm}
  Let $h$ be the highest common factor of the rank $r \geq 1$ and the degree $d \in \integers$. There is a birationally linear map of stacks
  \begin{displaymath} \xymatrix{
    \mu: \Bunstack_{r, d} \ar@{-->}[r] & \Bunstack_{h, 0}
  } \end{displaymath}
  and an isomorphism between the Picard \emph{schemes} $\Pic^d(C)$ and $\Pic^0(C)$ such that the following diagram commutes:
  \begin{equation} \label{preserve_det} \xymatrix{
    \Bunstack_{r, d} \ar@{-->}[r]^{\mu} \ar[d]_{\det} & \Bunstack_{h, 0} \ar[d]^{\det}\\
    \Pic^d(C) \ar[r]^{\sim} & \Pic^0(C)
  } \end{equation}
\end{thm}

\begin{rem} \upshape
  One cannot expect an isomorphism of Picard \emph{stacks} here: If (\ref{preserve_det}) were a commutative diagram of stacks, then choosing a general vector
  bundle $E$ of rank $r$ and degree $d$ would yield a commutative diagram of automorphism groups
  \begin{displaymath} \xymatrix{
    \Gm \ar[r]^{\sim} \ar[d]_{(\_)^r} & \Gm \ar[d]^{(\_)^h}\\
    \Gm \ar[r]^{\sim} & \Gm
  } \end{displaymath}
  which is impossible for $r \neq h$.
\end{rem}

\begin{rem} \upshape \label{scalar_autos}
  In the theorem, we can furthermore achieve that $\mu$ preserves scalar automorphisms in the following sense:

  Let $E$ and $E' = \mu(E)$ be vector bundles over $C$ that correspond to a general point in $\Bunstack_{r, d}$ and its image in $\Bunstack_{h, 0}$. Then $E$ and $E'$ are stable
  (because we have assumed $g \geq 2$) and hence simple. The rational map $\mu$ induces a morphism of algebraic groups
  \begin{displaymath}
    \mu^E: \Gm = \bbAut( E) \longto \bbAut( E') = \Gm
  \end{displaymath}
  which is an isomorphism because $\mu$ is birationally linear. Thus $\mu^E$ is either the identity or $\lambda \mapsto \lambda^{-1}$; it is independent of $E$ because $\Bunstack_{r, d}$ is
  irreducible. Modifying $\mu$ by the automorphism $E' \mapsto E'^{\dual}$ of $\Bunstack_{h, 0}$ if necessary, we can achieve that $\mu^E$ is the identity for every general $E$.
\end{rem}

Clearly, the map $\mu$ in the theorem restricts to a birationally linear map between the dense open substacks of \emph{stable} vector bundles. But any rational (resp.
birational, resp. birationally linear) map between these induces a rational (resp. birational, resp. birationally linear) map between the corresponding coarse moduli
schemes; cf. proposition \ref{descent} in the appendix for details. Hence the original theorem of King and Schofield follows:
\begin{cor}[King--Schofield] \label{coarse}
  Let $\coarseM_{r, d}$ be the coarse moduli scheme of stable vector bundles of rank $r$ and degree $d$ on $C$. Then there is a birationally linear map
  of schemes
  \begin{displaymath} \xymatrix{
    \mu: \coarseM_{r, d} \ar@{-->}[r] & \coarseM_{h, 0}.
  } \end{displaymath}
\end{cor}
Of course, this is just a reformulation of the theorem \ref{introthm} mentioned in the introduction.
\begin{rem} \upshape
  As mentioned before, King and Schofield also prove that the rational map $\xymatrix{ \mu: \coarseM_{r, d} \ar@{-->}[r] & \coarseM_{h, 0} }$ preserves their Brauer class $\psi_{r, d}$. This
  is equivalent to the condition that $\mu$ induces a rational map between the corresponding $\Gm$-gerbes, i.\,e. a rational map $\xymatrix{ \Bunstack_{r, d} \ar@{-->}[r] & \Bunstack_{h, 0} }$
  that preserves scalar automorphisms in the sense of remark \ref{scalar_autos}.
\end{rem}

We recall the consequences concerning the rationality of $\coarseM_{r, L}$. Because the diagram (\ref{preserve_det}) commutes, $\mu$ restricts to a rational map between fixed determinant moduli
schemes; thus one obtains:
\begin{cor}[King--Schofield] \label{rational}
  Let $L$ be a line bundle on $C$, and let $\coarseM_{r, L}$ be the coarse moduli scheme of stable vector bundles of rank $r$ and determinant $L$ on $C$.
  Then there is a birationally linear map of schemes
  \begin{displaymath} \xymatrix{
    \mu: \coarseM_{r, L} \ar@{-->}[r] & \coarseM_{h, \O}
  } \end{displaymath}
  where $h$ is the highest common factor of $r$ and $\deg(L)$.
\end{cor}
In particular, $\coarseM_{r, L}$ is rational if the rank $r$ and the degree $\deg(L)$ are coprime; this proves the conjecture mentioned in the introduction. More generally,
it follows that $\coarseM_{r, L}$ is rational if $\coarseM_{h, \O}$ is. For $h \geq 2$, it seems to be still an open question whether $\coarseM_{h, \O}$ is rational or not.

\section{Grassmannian bundles} \label{grassmannians}
Let $\stackV$ be a vector bundle over a dense open substack $\stackU \subseteq \Bunstack_{r, d}$. Recall that a part of this datum is a functor from the groupoid
$\stackU(k)$ to the groupoid of vector spaces over $k$. So for each appropriate vector bundle $E$ over $C$, we do not only get a vector space $\stackV_E$ over $k$,
but also a group homomorphism $\Aut_{\O_C}(E) \to \Aut_k( \stackV_E)$. Note that both groups contain the scalars $k^*$.
\begin{defn} \label{weight}
  A vector bundle $\stackV$ over a dense open substack $\stackU \subseteq \Bunstack_{r, d}$ \emph{has weight $w \in \integers$} if the diagram
  \begin{displaymath} \xymatrix{
    k^* \ar@{^{(}->}[r] \ar[d]^{(\_)^w} & \Aut_{\O_C}(E) \ar[d]\\
    k^* \ar@{^{(}->}[r]                 & \Aut_k( \stackV_E)
  } \end{displaymath}
  commutes for all vector bundles $E$ over $C$ that are objects of the groupoid $\stackU(k)$. 
\end{defn}

\begin{example}
  The trivial vector bundle $\O^n$ over $\Bunstack_{r, d}$ has weight $0$.
\end{example}

We denote by $\Euniv$ the universal vector bundle over $C \times \Bunstack_{r, d}$, and by $\Euniv_p$ its restriction to $\{p\} \times \Bunstack_{r, d}$ for some point
$p \in C(k)$.
\begin{example}
  $\Euniv_p$ is a vector bundle of weight $1$ on $\Bunstack_{r, d}$, and its dual $(\Euniv_p)^{\dual}$ is a vector bundle of weight $-1$.
\end{example}

For another example, we fix a vector bundle $F$ over $C$. By semicontinuity, there is an open substack $\stackU \subseteq \Bunstack_{r, d}$ that
parameterizes vector bundles $E$ of rank $r$ and degree $d$ over $C$ with $\Ext^1( F, E) = 0$; we assume $\stackU \neq \emptyset$. The vector spaces
$\Hom( F, E)$ are the fibres of a vector bundle $\Hom( F, \Euniv)$ over $\stackU$ according to Grothendieck's theory of cohomology and base change in EGA III.

Similarly, there is a vector bundle $\Hom( \Euniv, F)$ defined over an open substack of $\Bunstack_{r, d}$ whose fibre over
any point $[E]$ with $\Ext^1( E, F) = 0$ is the vector space $\Hom( E, F)$.
\begin{example}
  $\Hom( F, \Euniv)$ is a vector bundle of weight $1$, and $\Hom( \Euniv, F)$ is a vector bundle of weight $-1$.
\end{example}

Note that any vector bundle of weight $0$ over an open substack $\stackU \subseteq \Bunstack_{r, d}$ contained in the stable locus descends to a vector
bundle over the corresponding open subscheme $\coarseU \subseteq \coarseM_{r, d}$ of the coarse moduli scheme, cf. proposition \ref{descent}. Vector
bundles of nonzero weight do not descend to the coarse moduli scheme.

\begin{prop} \label{decompose}
  Consider all vector bundles $\stackV$ of fixed weight $w$ over dense open substacks of a fixed stack $\Bunstack_{r, d}$. Assume that $\stackV_0$ has minimal rank among
  them. Then every such $\stackV$ is generically isomorphic to $\stackV_0^n$ for some $n$.
\end{prop}
\begin{proof}
  The homomorphism bundles $\Endbdl( \stackV_0)$ and $\Hombdl( \stackV_0, \stackV)$ are vector bundles of weight $0$ over dense open substacks of $\Bunstack_{r, d}$. Hence
  they descend to vector bundles $A$ and $M$ over dense open subschemes of $\coarseM_{r, d}$, cf. proposition \ref{descent}. The algebra structure on $\Endbdl( \stackV_0)$
  and its right(!) action on $\Hombdl( \stackV_0, \stackV)$ also descend; they turn $A$ into an Azumaya algebra and $M$ into a right $A$-module.

  In particular, the generic fibre $M_K$ is a right module under the central simple algebra $A_K$ over the function field $K := k( \coarseM_{r, d})$. By our choice of
  $\stackV_0$, there are no nontrivial idempotent elements in $A_K$; hence $A_K$ is a skew field.
  
  We have just constructed a functor $\stackV \mapsto M_K$ from the category in question to the category of finite-dimensional right vector spaces over $A_K$. This functor
  is a Morita equivalence; its inverse is defined as follows:

  Given such a right vector space $M_K$ over $A_K$, we can extend it to a right $A$-module $M$ over a dense open subscheme of $\coarseM_{r, d}$, i.\,e. to a right
  $\Endbdl( \stackV_0)$-module of weight $0$ over a dense open substack of $\Bunstack_{r, d}$; we send $M_K$ to the vector bundle of weight $w$
  \begin{displaymath}
    \stackV := M \otimes_{\Endbdl( \stackV_0)} \stackV_0.
  \end{displaymath}

  Using this Morita equivalence, the proposition follows from the corresponding statement for right vector spaces over $A_K$.
\end{proof}
\begin{cor} \label{rank_h}
  There is a vector bundle of weight $w = 1$ (resp. $w = -1$) and rank $h = \hcf(r, d)$ over a dense open substack of $\Bunstack_{r, d}$.
\end{cor}
\begin{proof}
  Because the case of weight $w = -1$ follows by dualizing the vector bundles, we only consider vector bundles of weight $w = 1$. Here $\Euniv_p$ is a vector bundle of
  rank $r$ over $\Bunstack_{r, d}$, and $\Hom( L^{\dual}, \Euniv)$ is a vector bundle of rank $r( 1 - g + \deg(L)) + d$ over a dense open substack if $L$ is a sufficiently
  ample line bundle on $C$. Consequently, the rank of $\stackV_0$ divides $r$ and $r( 1 - g + \deg(L)) + d$; hence it also divides their highest common factor $h$.
\end{proof}

To each vector bundle $\stackV$ over a dense open substack $\stackU \subseteq \Bunstack_{r, d}$, we can associate a Grassmannian bundle
\begin{displaymath}
  \Gr_j( \stackV) \longto \stackU \subseteq \Bunstack_{r, d}.
\end{displaymath}
By definition, $\Gr_j( \stackV)$ is the moduli stack of those vector bundles $E$ over $C$ which are parameterized by $\stackU$, endowed with a
$j$-dimensional vector subspace of $\stackV_E$. $\Gr_j( \stackV)$ is again a smooth Artin stack locally of finite type over $k$, and
its canonical morphism to $\stackU$ is representable by Grassmannian bundles of schemes.

If $\stackV$ is a vector bundle of some weight, then all scalar automorphisms of $E$ preserve all vector subspaces of $\stackV_E$. This means that the automorphism groups
of the groupoid $\Gr_j( \stackV)(k)$ also contain the scalars $k^*$. In particular, it makes sense to say that a vector bundle over $\Gr_j( \stackV)$ has weight
$w \in \integers$: There is an obvious way to generalize definition \ref{weight} to this situation.

To give some examples, we fix a point $p \in C(k)$. Let $\Parstack_{r, d}^m$ be the moduli stack of rank $r$, degree $d$ vector bundles $E$ over $C$ endowed with a
quasiparabolic structure of multiplicity $m$ over $p$. Recall that such a quasiparabolic structure is just a coherent subsheaf $E' \subseteq E$ with the property that
$E/E'$ is isomorphic to the skyscraper sheaf $\O_p^m$.
\begin{example} \label{hecke_1}
  $\Parstack_{r, d}^m$ is canonically isomorphic to the Grassmannian bundle $\Gr_m( (\Euniv_p)^{\dual})$ over $\Bunstack_{r, d}$.
\end{example}
Here we have regarded a quasiparabolic vector bundle $E^{\bullet} = (E' \subseteq E)$ as the vector bundle $E$ together with a dimension $m$ quotient of the
fibre $E_p$. But we can also regard it as the vector bundle $E'$ together with a dimension $m$ vector subspace in the fibre at $p$ of the twisted vector bundle $E'(p)$.
Choosing a trivialization of the line bundle $\O_C(p)$ over $p$, we can identify the fibres of $E'(p)$ and $E'$ at $p$; hence we also obtain:
\begin{example} \label{hecke_2}
  $\Parstack_{r, d}^m$ is isomorphic to the Grassmannian bundle $\Gr_m(\Eprimeuniv_p)$ over $\Bunstack_{r, d-m}$ where $\Eprimeuniv$ is the universal vector bundle over
  $C \times \Bunstack_{r, d-m}$.
\end{example}

These two Grassmannian bundles
\begin{displaymath}
  \Bunstack_{r, d} \otgnol[\theta_1] \Parstack_{r, d}^m \longto[\theta_2] \Bunstack_{r, d-m}
\end{displaymath} 
form a correspondence between $\Bunstack_{r, d}$ and $\Bunstack_{r, d-m}$, the \emph{Hecke correspondence}. Its effect on the determinant line bundles is given by
\begin{equation} \label{hecke_det}
  \det \theta_1( E^{\bullet}) = \det( E) \cong \O_C(mp) \otimes \det( E') = \O_C(mp) \otimes \det \theta_2( E^{\bullet})
\end{equation}
for each parabolic vector bundle $E^{\bullet} = (E' \subseteq E)$ with multiplicity $m$ at $p$.

\begin{prop} \label{map}
  Let $\stackV$ and $\stackW$ be two vector bundles of the same weight $w$ over dense open substacks of $\Bunstack_{r, d}$. If $j \leq \rank(\stackW)
  \leq \rank(\stackV)$, then there is a birationally linear map
  \begin{displaymath} \xymatrix{
    \rho: \Gr_j( \stackV) \ar@{-->}[r] & \Gr_j( \stackW)
  } \end{displaymath}
  over $\Bunstack_{r, d}$.
\end{prop}
\begin{proof}
  According to proposition \ref{decompose}, there is a vector bundle $\stackW'$ of weight $w$ such that $\stackV \cong \stackW \oplus \stackW'$ over some
  dense open substack $\stackU \subseteq \Bunstack_{r, d}$. We may assume without loss of generality that $\stackU$ is contained in the stable locus
  and denote by $\coarseU \subseteq \coarseM_{r, d}$ the corresponding open subscheme, cf. proposition \ref{descent}.

  We use the following simple fact from linear algebra: If $W$ and $W'$ are vector spaces over $k$ with $\dim(W) \geq j$, then every $j$-dimensional
  vector subspace of $W \oplus W'$ whose image $S$ in $W$ also has dimension $j$ is the graph of a unique linear map $S \to W'$.

  This means that $\Gr_j( W \oplus W')$ contains as a dense open subscheme the total space of the vector bundle $\Hom( S^{\univ}, W')$ over $\Gr_j( W)$ where $S^{\univ}$
  is the universal subbundle of the constant vector bundle $W$ over $\Gr_j( W)$.

  In our stack situation, these considerations imply that $\Gr_j( \stackV)$ is birational to the total space of the vector bundle $\Hombdl( \Suniv, \stackW')$ over
  $\Gr_j( \stackW)$ where $\Suniv$ is the universal subbundle of the pullback of $\stackW$ over $\Gr_j( \stackW)$. This defines the rational map $\rho$.

  The vector bundle $\Hombdl( \Suniv, \stackW')$ has weight $0$ because $\Suniv$ and $\stackW'$ both have weight $w$. Since the scalars act trivially,
  we can descend $\Gr_j( \stackW)$ and this vector bundle over it to a Grassmannian bundle over $\coarseU$ and a vector bundle over it, cf. proposition
  \ref{descent}. In particular, our homomorphism bundle is trivial over a dense open substack of $\Gr_j( \stackW)$. This proves that $\rho$ is birationally linear.
\end{proof}
\begin{cor} \label{gr_rational}
  Let $\stackV$ be a vector bundle of weight $w = \pm1$ over a dense open substack of $\Bunstack_{r, d}$. If $j$ is divisible by $\hcf(r, d)$, then the Grassmannian bundle
  \begin{displaymath}
    \Gr_j( \stackV) \longto \Bunstack_{r, d}
  \end{displaymath}
  is birationally linear.
\end{cor}
\begin{proof}
  By corollary \ref{rank_h}, there is a vector bundle $\stackW$ of weight $w$ and rank $j$. Due to the proposition, $\Gr_j( \stackV)$ is birationally linear over
  $\Gr_j( \stackW) \simeq \Bunstack_{r, d}$.
\end{proof}

\section{Proof of theorem \ref{mainthm}} \label{proof}

The aim of this section is to prove theorem \ref{mainthm}, i.\,e. to construct the birationally linear map $\xymatrix{\mu: \Bunstack_{r, d}
\ar@{-->}[r] & \Bunstack_{h, 0}}$ where $h$ is the highest common factor of the rank $r$ and the degree $d$. We proceed by induction on $r/h$.

For $r = h$, the theorem is trivial: Tensoring with an appropriate line bundle defines even an isomorphism of stacks
$\xymatrix{\Bunstack_{r, d} \ar[r]^{\sim} & \Bunstack_{h, 0}}$ with the required properties. Thus we may assume $r > h$.
\begin{lemma}
  There are unique integers $r_F$ and $d_F$ that satisfy
  \begin{equation} \label{eulerchar}
    (1 - g) r_F r + r_F d - r d_F = h
  \end{equation}
  and 
  \begin{equation} \label{smallerrank}
    r < hr_F < 2r.
  \end{equation}
\end{lemma}
\begin{proof}
  (\ref{eulerchar}) has an integer solution $r_F, d_F$ because $h$ is also the highest common factor of $r$ and $(1-g)r + d$; here $r_F$ is unique modulo
  $r/h$. Furthermore, $r_F$ is nonzero modulo $r/h$ since $- r d_F = h$ has no solution. Hence precisely one of the solutions $r_F, d_F$ of
  (\ref{eulerchar}) also satisfies (\ref{smallerrank}).
\end{proof}

We fix $r_F$, $d_F$ and define
\begin{displaymath}
  r_1 := h r_F - r, \qquad d_1 := h d_F - d, \qquad h_1 := \mathrm{hcf}( r_1, d_1).
\end{displaymath}
Then $r_1 < r$, and $h_1$ is a multiple of $h$. In particular, $r_1/h_1 < r/h$.

\begin{lemma} \label{F-generic}
  There is an exact sequence
  \begin{equation} \label{F-sequence}
    0 \longto E_1 \longto F \otimes_k V \longto E \longto 0
  \end{equation}
  where $E_1$, $F$, $E$ are vector bundles over $C$ and $V$ is a vector space over $k$ with
  \begin{displaymath} \begin{array}{r@{\qquad}r@{\qquad}r@{\qquad}r}
    \rank(E_1) = r_1, & \rank(F) = r_F, & \rank(E) = r,  & \dim(V) = h\\
    \deg (E_1) = d_1, & \deg (F) = d_F, & \deg (E) = d\, &            \\
  \end{array} \end{displaymath}
  such that the following two conditions are satisfied:
  \begin{itemize}
   \item[i)] $\Ext^1( F, E) = 0$, and the induced map $V \to \Hom( F, E)$ is bijective.
   \item[ii)] $\Ext^1( E_1, F) = 0$, and the induced map $V^{\dual} \to \Hom( E_1, F)$ is injective.
  \end{itemize}
\end{lemma}
\begin{proof}
  We may assume $h = 1$ without loss of generality: If there is such a sequence for $r/h$ and $d/h$ instead of $r$ and $d$, then the direct sum of $h$
  copies is the required sequence for $r$ and $d$.

  By our choice of $r_F$ and $d_F$ and Riemann-Roch, all vector bundles $F$ and $E$ of these ranks and degrees satisfy
  \begin{displaymath}
    \chi( F, E) := \dim_k \Hom( F, E) - \dim_k \Ext^1( F, E) = h = 1.
  \end{displaymath}
  If $F$ and $E$ are general, then
  \begin{displaymath}
    \Hom( F, E) \cong k \qquad\text{and}\qquad \Ext^1( F, E) = 0
  \end{displaymath}
  according to a theorem of Hirschowitz \cite[section 4.6]{hirschowitz}, and there is a surjective map $F \to E$ by an argument of Russo and Teixidor
  \cite{russo-teixidor}. Thus we obtain an exact sequence
  \begin{equation} \label{F-sequence2}
    0 \longto E_1 \longto F \longto E \longto 0
  \end{equation}
  that satisfies condition i (with $V = k$).

  (For the convenience of the reader, a proof of the cited results is given in the appendix, cf. theorem \ref{nonspecial}.)

  Furthermore, all vector bundles of the given ranks and degrees satisfy
  \begin{displaymath}
    \chi( E_1, F) = \chi( F, E) - \chi( E, E) + \chi( E_1, E_1) > \chi( F, E) = h = 1
  \end{displaymath}
  because $r_1 < r$. Now we can argue as above: For general $E_1$ and $F$, we have $\Ext^1( E_1, F) = 0$ by Hirschowitz, and there is an injective map
  $E_1 \to F$ with torsionfree cokernel by Russo-Teixidor; cf. also theorem \ref{nonspecial} in the appendix. Thus we obtain an exact sequence
  (\ref{F-sequence2}) that satisfies condition ii (with $V = k$).

  Finally, we consider the moduli stack of all exact sequences (\ref{F-sequence2}) of vector bundles with the given ranks and degrees. As explained
  in the appendix (cf. corollary \ref{integral}), it is an \emph{irreducible} algebraic stack locally of finite type over $k$. But i and ii are open
  conditions, so there is a sequence that satisfies both.
\end{proof}

From now on, let $F$ be a fixed vector bundle of rank $r_F$ and degree $d_F$ that occurs in such an exact sequence (\ref{F-sequence}).
\begin{defn}
  The rational map of stacks
  \begin{displaymath} \xymatrix{
    \lambda_F: \Bunstack_{r, d} \ar@{-->}[r] & \Bunstack_{r_1, d_1}
  } \end{displaymath}
  is defined by sending a general rank $r$, degree $d$ vector bundle $E$ over $C$ to the kernel of the natural evaluation map
  \begin{displaymath}
    \epsilon_E: \Hom(F, E) \otimes_k F \longto E.
  \end{displaymath}
\end{defn}
We check that this does define a rational map. Let $\stackU_F \subseteq \Bunstack_{r, d}$ be the open substack that parameterizes all $E$
for which $\Ext^1( F, E) = 0$ and $\epsilon_E$ is surjective. Then the $\epsilon_E$ are the restrictions of a surjective morphism $\epsilon^{\univ}$
of vector bundles over $C \times \stackU_F$. So the kernel of $\epsilon^{\univ}$ is also a vector bundle; it defines a morphism
$\lambda_F: \stackU_F \to \Bunstack_{r_1, d_1}$. This gives the required rational map because $\stackU_F$ is nonempty by our choice of $F$.

For later use, we record the effect of $\lambda_F$ on determinant line bundles:
\begin{equation} \label{lambda_det}
  \det \lambda_F( E) \cong \det( F)^{\otimes h} \otimes \det( E)^{\dual}.
\end{equation}

Following \cite{king-schofield}, the next step is to understand the fibres of $\lambda_F$. We denote by $\Hom( \Euniv_1, F)$ the vector bundle over an open substack of
$\Bunstack_{r_1, d_1}$ whose fibre over any point $[E_1]$ with $\Ext^1( E_1, F) = 0$ is the vector space $\Hom( E_1, F)$.
\begin{prop} \label{is_grass}
  $\Bunstack_{r, d}$ is over $\Bunstack_{r_1, d_1}$ naturally birational to the Grassmannian bundle $\Gr_h( \Hom( \Euniv_1, F))$.
\end{prop}
\begin{proof}
  If $E$ is a rank $r$, degree $d$ vector bundle over $C$ for which $\Ext^1( F, E) = 0$ and the above map $\epsilon := \epsilon_E$ is surjective, then
  the exact sequence
  \begin{displaymath}
    0 \longto \ker( \epsilon) \longto \Hom( F, E) \otimes _k F \longto[\epsilon] E \longto 0
  \end{displaymath}
  satisfies the condition i of the previous lemma. This identifies the above open substack $\stackU_F \subseteq \Bunstack_{r, d}$ with the moduli stack
  of all exact sequences (\ref{F-sequence}) that satisfy i.

  Similarly, let $\stackU_F' \subseteq \Gr_h( \Hom( \Euniv_1, F))$ be the open substack that parameterizes all pairs $(E_1, S \subseteq \Hom( E_1, F))$
  for which $\Ext^1( E_1, F) = 0$ and the natural map $\alpha: E_1 \to S^{\dual} \otimes_k F$ is injective with torsionfree cokernel. For such a
  pair $(E_1, S)$, the exact sequence
  \begin{displaymath}
    0 \longto E_1 \longto[\alpha] S^{\dual} \otimes _k F \longto \coker( \alpha) \longto 0
  \end{displaymath}
  satisfies the condition ii of the previous lemma. This identifies $\stackU_F'$ with the moduli stack of all exact sequences (\ref{F-sequence}) that satisfy ii.

  Hence both $\Bunstack_{r, d}$ and $\Gr_h( \Hom( \Euniv_1, F))$ contain as an open substack the moduli stack $\stackU_F''$ of all exact sequences
  (\ref{F-sequence}) that satisfy both conditions i and ii. But $\stackU_F''$ is nonempty by our choice of $F$, so it is dense in both stacks; thus
  they are birational over $\Bunstack_{r_1, d_1}$.
\end{proof}

Still following \cite{king-schofield}, the proof of theorem \ref{mainthm} can now be summarized in the following diagram; it is explained below.
\begin{displaymath} \xymatrix{
  \Bunstack_{r, d}      \ar@{-->}[dr]_{\lambda_F} \ar@{-->}[r]^{\rho} & \Gr_h( \stackW) \ar[d] \ar@{-->}^{\tilde{\mu_1}}[r] &
    \Parstack_{h_1,0}^h \ar[d]^{\theta_1}         \ar[r]^-{\theta_2} & \Bunstack_{h_1, -h} \ar@{-->}[r]^{\mu_2} & \Bunstack_{h, 0}\\
  & \Bunstack_{r_1, d_1} \ar@{-->}[r]_{\mu_1} & \Bunstack_{h_1, 0}
} \end{displaymath}
Here $\mu_1$ and $\mu_2$ are the birationally linear maps given by the induction hypothesis. $(\theta_1, \theta_2)$ is the Hecke correspondence explained in the previous
section; note that $\theta_2$ is birationally linear by corollary \ref{gr_rational}.

The square in this diagram is cartesian, so $\tilde{\mu_1}$ is the pullback of $\mu_1$, and $\stackW := \mu_1^* (\Euniv_p)^{\dual}$ is the pullback of the vector bundle
$(\Euniv_p)^{\dual}$ over $\Bunstack_{h_1, 0}$ to which $\theta_1$ is the associated Grassmannian bundle. Using remark \ref{scalar_autos}, we may assume that $\mu_1$ preserves
scalar automorphisms, i.\,e. that $\stackW$ has the same weight $-1$ as $(\Euniv_p)^{\dual}$. Then we can apply proposition \ref{map} to obtain the birationally linear map $\rho$.
Now we have the required birationally linear map
\begin{displaymath} \xymatrix{
  \mu := \mu_2 \circ \theta_2 \circ \tilde{\mu_1} \circ \rho: \Bunstack_{r, d} \ar@{-->}[r] & \Bunstack_{h, 0};
} \end{displaymath}
it satisfies the determinant condition in theorem \ref{mainthm} due to equations (\ref{lambda_det}), (\ref{hecke_det}) and the corresponding induction hypothesis on
$\mu_1, \mu_2$.

\begin{appendix}
\section{Moduli stacks of sheaves on curves} \label{stacks}
This section summarizes some well-known basic properties of moduli stacks of vector bundles and more generally coherent sheaves on curves. For the general theory of
algebraic stacks, we refer the reader to \cite{laumon} or the appendix of \cite{vistoli}. We prove that the moduli stacks in question are algebraic, smooth and irreducible.
Then we discuss descent to the coarse moduli scheme. Finally, we deduce Hirschowitz' theorem \cite{hirschowitz} and a refinement by Russo and Teixidor \cite{russo-teixidor}
about morphisms between general vector bundles. 

Recall that we have fixed an algebraically closed field $k$ and a connected smooth projective curve $C/k$ of genus $g$. We say that a coherent
sheaf $F$ on $C$ has \emph{type} $t = (r, d)$ if its rank $\rank(F)$ (at the generic point of $C$) equals $r$ and its degree $\deg(F)$ equals $d$.

If $F'$ and $F$ are coherent sheaves of types $t = (r, d)$ and $t' = (r', d')$ on $C$, then the Euler characteristic
\begin{displaymath}
  \chi(F', F) := \dim_k \Hom( F', F) - \dim_k \Ext^1( F', F)
\end{displaymath}
satisfies the Riemann-Roch theorem $\chi( F', F)$ = $\chi( t', t)$ with
\begin{displaymath}
  \chi(t', t) := (1- g) r' r + r' d - r d'.
\end{displaymath}
Note that $Ext^n( F', F) $ vanishes for all $n \geq 2$ since $\dim( C) = 1$.

We denote by $\Cohstack_t$ the moduli stack of coherent sheaves $F$ of type $t$ on $C$. More precisely, $\Cohstack_t( S)$ is for each $k$-scheme $S$
the groupoid of all coherent sheaves on $C \times S$ which are flat over $S$ and whose fibre over every point of $S$ has type $t$.

Now assume $t = t_1 + t_2$. We denote by $\Extstack( t_2, t_1)$ the moduli stack of exact sequences of coherent sheaves on $C$
\begin{displaymath}
  0 \to F_1 \to F \to F_2 \to 0
\end{displaymath}
where $F_i$ has type $t_i = (r_i, d_i)$ for $i = 1, 2$. This means that $\Extstack( t_2, t_1)(S)$ is for each $k$-scheme $S$
the groupoid of short exact sequences of coherent sheaves on $C \times S$ which are flat over $S$ and fibrewise of the given types.
\begin{prop} \label{algebraic}
  The stacks $\Cohstack_t$ and $\Extstack( t_2, t_1)$ are algebraic in the sense of Artin and locally of finite type over $k$.
\end{prop}
\begin{proof}
  Let $\O(1)$ be an ample line bundle on $C$. For $n \in \integers$, we denote by
  \begin{displaymath}
    \Cohstack^n_t \subseteq \Cohstack_t \qquad (\text{resp. } \Extstack( t_2, t_1)^n \subseteq  \Extstack( t_2, t_1))
  \end{displaymath}
  the open substack that parameterizes coherent sheaves $F$ on $C$ (resp. exact sequences $0 \to F_1 \to F \to F_2 \to 0$ of coherent sheaves on $C$) such that the twist
  $F(n) = F \otimes \O(1)^{\otimes n}$ is generated by global sections and $\cohom^1( F(n)) = 0$.

  By Grothendieck's theory of Quot-schemes, there is a scheme $\Quot^n_t$ of finite type over $k$ that parameterizes such coherent sheaves $F$ together
  with a basis of the $k$-vector space $\cohom^0( F( n))$. Moreover, there is a relative Quot-scheme $\Flag( t_2, t_1)^n$ of finite type over $\Quot^n_t$ that
  parameterizes such exact sequences $0 \to F_1 \to F \to F_2 \to 0$ together with a basis of $\cohom^0( F( n))$.

  Let $N$ denote the dimension of $\cohom^0( F( n))$. According to Riemann-Roch, $N$ depends only on $t$, $n$ and the ample line bundle $\O(1)$, but not on $F$.

  Changing the chosen basis defines an action of $\GL( N)$ on $\Quot^n_t$, and $\Cohstack^n_t$ is precisely the stack quotient $\Quot^n_t/\GL( N)$. Similarly,
  $\Extstack( t_2, t_1)^n$ is precisely the stack quotient $\Flag( t_2, t_1)^n/\GL( N)$. Hence these two stacks are algebraic and of finite type over $k$.

  By the ampleness of $\O(1)$, the $\Cohstack^n_t$ (resp. $\Extstack( t_2, t_1)^n$) form an open covering of $\Cohstack_t$ (resp. $\Extstack( t_2, t_1)$).
\end{proof}

\begin{rem} \upshape
  In general, $\Cohstack_t$ is not quasi-compact because the family of all coherent sheaves on $C$ of type $t$ is not bounded.
\end{rem}

\begin{prop} \label{smooth}
  \begin{itemize}
   \item[i)] $\Cohstack_t$ is smooth of dimension $-\chi(t, t)$ over $k$.
   \item[ii)] If we assign to each exact sequence $0 \to F_1 \to F \to F_2 \to 0$ the two sheaves $F_1, F_2$, then the resulting morphism of algebraic stacks
    \begin{displaymath}
      \Extstack( t_2, t_1) \longto \Cohstack( t_1) \times \Cohstack( t_2)
    \end{displaymath}
    is smooth of relative dimension $-\chi( t_2, t_1)$, and all its fibres are connected.
   \item[iii)] $\Extstack( t_2, t_1)$ is smooth of dimension $-\chi(t_2, t_2) - \chi( t_2, t_1) - \chi( t_1, t_1)$ over $k$.
  \end{itemize}
\end{prop}
\begin{proof}
  i) By standard deformation theory, $\Hom( F, F)$ is the automorphism group of any infinitesimal deformation of the coherent sheaf $F$,
  $\Ext^1( F, F)$ classifies such deformations, and $\Ext^2( F, F)$ contains the obstructions against extending deformations infinitesimally,
  cf. \cite[2.A.6]{huybrechts-lehn}. Because $\Ext^2$ vanishes, deformations of $F$ are unobstructed and hence $\Cohstack_t$ is smooth;
  its dimension at $F$ is then $\dim \Ext^1( F, F) - \dim \Hom( F, F) = -\chi( t, t)$.

  ii) The fibre of this morphism over $[F_1, F_2]$ is the moduli stack of all extensions of $F_2$ by $F_1$; it is the stack quotient of the affine space $\Ext^1( F_2, F_1)$
  modulo the trivial action of the algebraic group $\Hom( F_2, F_1)$. Hence this fibre is smooth of dimension $-\chi( t_2, t_1)$ and connected.

  More generally, consider a scheme $S$ of finite type over $k$ and a morphism $S \to \Cohstack( t_1) \times \Cohstack( t_2)$; let $F_1$ and $F_2$ be the
  corresponding coherent sheaves over $C \times S$. By EGA III, the object $\RHom( F_2, F_1)$ in the derived category of coherent sheaves on $S$ can
  locally be represented by a complex of length one $V^0 \longto[\delta] V^1$ where $V^0, V^1$ are vector bundles. This means that the pullback of
  $\Extstack( t_2, t_1)$ to $S$ is locally the stack quotient of the total space of $V^1$ modulo the action of the algebraic group $V^0 \big/ S$ given
  by $\delta$. Hence this pullback is smooth over $S$; this proves ii.

  iii) follows from i and ii.
\end{proof}

\begin{prop} \label{connected}
  The stacks $\Cohstack_t$ and $\Extstack( t_2, t_1)$ are connected.
\end{prop}
\begin{proof}
  Proposition \ref{smooth} implies that $\Extstack( t_2, t_1)$ is connected if $\Cohstack_{t_1}$ and $\Cohstack_{t_2}$ are.
  We prove the connectedness of the latter by induction on the rank (and on the degree for rank zero).

  $\Cohstack_t$ is connected for $t = (0, 1)$ because there is a canonical surjection $C \to \Cohstack_t$; it sends a point $P$ to the skyscraper sheaf
  $\O_P$. Now consider $t = (0, d)$ with $d \geq 2$ and write $t = t_1 + t_2$. By induction hypothesis and \ref{smooth}, $\Extstack( t_1, t_2)$ is
  connected. But there is a canonical surjection $\Extstack( t_1, t_2) \to \Cohstack_t$; it sends an extension $0 \to F_1 \to F \to F_2 \to 0$ to the
  sheaf $F$. This shows that $\Cohstack_t$ is also connected; now we have proved all connectedness assertions in rank zero.

  If $F$ and $F'$ are two coherent sheaves on $C$ of type $t = (r, d)$ with $r \geq 1$, then there is a line bundle $L$ on $C$ such that both
  $L^{\dual} \otimes F$ and $L^{\dual} \otimes F'$ have a generically nonzero section. In other words, there are injective morphisms $L \hookrightarrow F$
  and $L \hookrightarrow F'$. Let $t_L$ be the type of $L$; then $F$ and $F'$ are both in the image of the canonical morphism $\Extstack( t - t_L, t_L)
  \to \Cohstack_t$. But $\Extstack( t - t_L, t_L)$ is connected by the induction hypothesis and \ref{smooth}. This shows that any two points $F$
  and $F'$ lie in the same connected component of $\Cohstack_t$, i.\,e. $\Cohstack_t$ is connected.
\end{proof}

\begin{cor} \label{integral}
  The stacks $\Cohstack_t$ and $\Extstack( t_2, t_1)$ are reduced and irreducible.
\end{cor}

The moduli stack $\Bunstack_t$ of vector bundles, the moduli stack $\Semistabstack_t$ of semi\-stable vector bundles and the moduli stack $\Stabstack_t$
of (geometrically) stable vector bundles on $C$ of type $t = (r, d)$ are open substacks
\begin{displaymath}
  \Stabstack_t \subseteq \Semistabstack_t \subseteq \Bunstack_t \subseteq \Cohstack_t.
\end{displaymath}
Hence these stacks are all irreducible and smooth of the same dimension $-\chi(t, t)$ if they are nonempty. $\Stabstack_t$ is known to be nonempty for
$g \geq 2$ and $r \geq 1$. Moreover, $\Semistabstack_t$ and $\Stabstack_t$ are quasi-compact (and thus of finite type) because the family of
(semi-)stable vector bundles of given type $t$ is bounded.

\begin{prop} \label{descent}
  Assume $g \geq 2$. Let $\Stabstack_t \longto \coarseM_t$ be the coarse moduli scheme of stable vector bundles of type $t$, and let $\stackV$ be a
  vector bundle of some weight $w \in \integers$ over an open substack $\stackU \subseteq \Stabstack_t$.
  \begin{itemize}
   \item[i)] $\stackU$ descends to an open subscheme $\coarseU \subseteq \coarseM_t$.
   \item[ii)] $\Gr_j( \stackV)$ descends to a (twisted) Grassmannian scheme $\coarseGr_j( \stackV)$ over $\coarseU$.
   \item[iii)] If $\stackV$ has weight $w = 0$, then it descends to a vector bundle over $\coarseU$.
   \item[iv)] More generally, any vector bundle of weight $0$ over $\Gr_j( \stackV)$ descends to a vector bundle over $\coarseGr_j( \stackV)$.
   \item[v)] Any birationally linear map of stacks $\xymatrix{\Stabstack_{t'} \ar@{-->}[r] & \Stabstack_t}$ induces a birationally linear
    map of schemes $\xymatrix{\coarseM_{t'} \ar@{-->}[r] & \coarseM_t}$.
  \end{itemize}
\end{prop}
\begin{proof}
  We continue to use the notation introduced in the proof of proposition \ref{algebraic}. By boundedness, there is an integer $n$ such that
  $\Stabstack_t$ is contained in $\Cohstack_t^n$; hence $\Stabstack_t = \Quot^{\stab}_t/\GL( N)$ where $\Quot^{\stab}_t \subseteq \Quot^n_t$ is
  the stable locus. Here the center of $\GL( N)$ acts trivially; by Geometric Invariant Theory \cite{git}, $\Quot^{\stab}_t$ is a principal
  $\PGL( N)$-bundle over $\coarseM_t$.

  i) Let $U \subseteq \Quot^{\stab}_t$ be the inverse image of $\stackU$. Then $U$ is a $\PGL( N)$-stable open subscheme in the total space of this
  principal bundle and hence the inverse image of a unique open subscheme $\coarseU \subseteq \coarseM_t$.

  ii) Let $V$ be the pullback of $\stackV$ to $U$; it is a vector bundle with $\GL( N)$-action. Hence its Grassmannian scheme $\Gr_j( V) \to U$ also
  carries an action of $\GL( N)$. But here the center acts trivially: $\lambda \cdot \id \in \GL(N)$ acts as the scalar $\lambda^w$ on the fibres of $V$
  and hence acts trivially on $\Gr_j( V)$. Thus we obtain an action of $\PGL( N)$ on our Grassmannian scheme $\Gr_j( V)$ over $U$. Since this action is
  free, $\Gr_j(V)$ descends to a Grassmannian bundle $\coarseGr_j( \stackV)$ over $\coarseU$ (which may be twisted, i.\,e. not Zariski-locally trivial).

  iii) The assumption $w=0$ means that the scalars in $\GL( N)$ act trivially on the fibres of $V$. Hence $\PGL( N)$ acts on $V$ over $U$ here. Again
  since this action is free, $V$ descends to a vector bundle over $\coarseU$.

  iv) Here weight $0$ means that the scalars in $\GL( N)$ act trivially on the pullback of the given vector bundle to $\Gr_j( V)$. Hence $\PGL( N)$ acts
  on this pullback; but it acts freely on the base $\Gr_j( V)$, so the vector bundle descends to $\coarseGr_j( \stackV)$.

  v) Such a birationally linear map can be represented by an isomorphism $\varphi: \stackU' \to \stackU$ between dense open substacks
  $\stackU' \subseteq \Stabstack_{t'} \times \mathbb{A}^?$ and $\stackU \subseteq \Stabstack_t$. By i, $\stackU$ corresponds to an open subscheme
  $\coarseU \subseteq \coarseM_t$; the proof of i shows that $\coarseU$ is a coarse moduli scheme for the stack $\stackU$. A straightforward
  generalization of this argument shows that $\stackU'$ corresponds to an open subscheme $\coarseU' \subseteq \coarseM_{t'} \times \mathbb{A}^?$ and
  that $\coarseU'$ is again a coarse moduli scheme for $\stackU'$. By the universal property of coarse moduli schemes, $\varphi$ induces an isomorphism
  $\coarseU' \to \coarseU$ and thus the required birationally linear map of schemes.
\end{proof}

\begin{thm}[Hirschowitz, Russo-Teixidor] \label{nonspecial}
  Assume $g \geq 2$. Let $F_1$ and $F_2$ be a general pair of vector bundles on $C$ with given types $t_1 = (r_1, d_1)$ and $t_2 = (r_2, d_2)$.
  \begin{itemize}
   \item[i)] If $\chi( t_1, t_2) \geq 0$, then $\dim \Hom( F_1, F_2) = \chi( t_1, t_2)$ and $\Ext^1( F_1, F_2) = 0$.
   \item[ii)] If $\chi( t_1, t_2) \geq 1$ and $r_1 > r_2$ (resp. $r_1 = r_2$, resp. $r_1 < r_2$), then a general morphism $F_1 \to F_2$ is
    surjective (resp. injective, resp. injective with torsionfree cokernel).
  \end{itemize}
\end{thm}
\begin{proof}
  The cases $r_1 = 0$ and $r_2 = 0$ are trivial, so we may assume $r_1, r_2 \geq 1$; then $\Stabstack_{t_1} \neq \emptyset \neq \Stabstack_{t_2}$.
  By semicontinuity, there is a dense open substack $\stackU \subseteq \Stabstack_{t_1} \times \Stabstack_{t_2}$ of stable vector bundles $F_1, F_2$
  with $\dim \Hom( F_1, F_2)$ minimal, say equal to $m$. According to Riemann-Roch, $m \geq \chi( t_1, t_2)$; part i of the theorem precisely claims
  that we have equality here.

  Let $\Hom( \Funiv_1, \Funiv_2)$ be the vector bundle of rank $m$ over $\stackU$ whose fibre over $F_1, F_2$ is $\Hom( F_1, F_2)$. By generic flatness
  (cf. EGA IV, \S6.9), there is a dense open substack $\stackV$ in the total space of $\Hom( \Funiv_1, \Funiv_2)$ such that the cokernel of the
  universal family of morphisms $F_1 \to F_2$ is flat over $\stackV$. Then its image and kernel are also flat over $\stackV$; we denote the types of
  cokernel, image and kernel by $t_Q = (r_Q, d_Q)$, $t = (r, d)$ and $t_K = (r_K, d_K)$.

  If $r = 0$, then the theorem clearly holds: In this case, the general morphism $\varphi: F_1 \to F_2$ has generic rank zero, so $\varphi = 0$; this
  means $m = 0$. Together with $m \geq \chi( t_1, t_2)$, this proves i and shows that the hypothesis of ii cannot hold here. Henceforth, we may thus
  assume $r \neq 0$.

  Note that $t_1 = t_K + t$ and $t_2 = t + t_Q$; moreover, we have a canonical morphism of moduli stacks
  \begin{displaymath}
    \Phi: \stackV \longto \Extstack( t, t_K) \times_{\Cohstack_t} \Extstack( t_Q, t)
  \end{displaymath}
  that sends a morphism $\varphi: F_1 \to F_2$ to the extensions
  \begin{displaymath}
    0 \to \ker( \varphi) \to F_1 \to \im( \varphi) \to 0 \quad\text{and}\quad 0 \to \im( \varphi) \to F_2 \to \coker( \varphi) \to 0.
  \end{displaymath}
  Conversely, two extensions $0 \to K \to F_1 \to I \to 0$ and $0 \to J \to F_2 \to Q \to 0$ together with an isomorphism $I \to J$ determine a
  morphism $\varphi: F_1 \to F_2$. Thus $\Phi$ is an isomorphism onto the open locus in $\Extstack \times_{\Cohstack} \Extstack$ where both extension
  sheaves $F_1, F_2$ are stable vector bundles and $\dim \Hom( F_1, F_2) = m$. Hence the stack dimensions coincide, i.\,e.
  \begin{displaymath}
    m - \chi( t_1, t_1) - \chi( t_2, t_2) = - \chi( t_1, t_K) - \chi( t, t) - \chi( t_Q, t_2).
  \end{displaymath}
  Since $\chi$ is biadditive, this is equivalent to
  \begin{equation} \label{excess}
    m - \chi( t_1, t_2) = - \chi( t_K, t_Q).
  \end{equation}
  In particular, $\chi( t_K, t_Q) \leq 0$ follows.

  Now suppose that $t_K$ and $t_Q$ were both nonzero. Since the general vector bundles $F_1$ and $F_2$ are stable, we then have
  \begin{displaymath}
    \frac{d_K}{r_K} < \frac{d_1}{r_1} < \frac{d}{r} < \frac{d_2}{r_2} < \frac{d_Q}{r_Q}.
  \end{displaymath}
  Using the assumption $\chi( t_1, t_2) \geq 0$, we get
  \begin{displaymath}
    \frac{\chi( t_K, t_Q)}{r_K r_Q} = 1 - g - \frac{d_K}{r_K} + \frac{d_Q}{r_Q} > 1 - g - \frac{d_1}{r_1} + \frac{d_2}{r_2} =
    \frac{\chi( t_1, t_2)}{r_1 r_2} \geq 0
  \end{displaymath}
  and hence $\chi( t_K, t_Q) > 0$. This contradiction proves $t_K = 0$ or $t_Q = 0$.

  (In some sense, this argument also covers the cases $r_K = 0$ and $r_Q = 0$. More precisely, $r_K = 0$ implies $t_K = 0$ because every rank zero
  coherent subsheaf of a vector bundle $F_1$ is trivial. On the other hand, $r_K \neq 0$ and $t_Q = (0, d_Q) \neq 0$ would imply
  $\chi( t_K, t_Q) = r_K d_Q > 0$ which is again a contradiction.)

  In particular, we get $\chi( t_K, t_Q) = 0$; together with equation (\ref{excess}), this proves part i of the theorem.

  If $r_1 > r_2$ (resp. $r_1 \leq r_2$), then $r_K > r_Q$ (resp. $r_K \leq r_Q$) and hence $r_K \neq 0 = r_Q$ (resp. $r_K = 0$); we have just seen 
  that this implies $t_Q = 0$ (resp. $t_K = 0$), i.\,e. the general morphism $\varphi: F_1 \to F_2$ is surjective (resp. injective).

  Furthermore, the morphism of stacks $\stackV \to \Cohstack_{t_Q}$ that sends a morphism $\varphi: F_1 \to F_2$ to its cokernel is smooth (due to
  the open embedding $\Phi$ and proposition \ref{smooth}). If $r_1 < r_2$, then $r_Q \geq 1$, so $\Bunstack_{t_Q}$
  is open and dense in $\Cohstack_{t_Q}$; this implies that the inverse image of $\Bunstack_{t_Q}$ is open and dense in $\stackV$, i.\,e. the general
  morphism $\varphi: F_1 \to F_2$ has torsionfree cokernel.
\end{proof}

\end{appendix}


\begin{thebibliography}{10}

\bibitem{ega}
A.~Grothendieck.
\newblock \'{E}l\'ements de g\'eom\'etrie alg\'ebrique ({EGA}).
\newblock {\em Inst. Hautes \'Etudes Sci. Publ. Math.}, 4, 8, 11, 17, 20, 24,
  28, 32, 1960--1967.

\bibitem{hirschowitz}
A.~Hirschowitz.
\newblock {P}robl\`emes de {B}rill-{N}oether en rang sup\`erieur.
\newblock \\ \verb|http://math.unice.fr/~ah/math/Brill/|.

\bibitem{huybrechts-lehn}
D.~Huybrechts and M.~Lehn.
\newblock {\em The geometry of moduli spaces of sheaves}.
\newblock Friedr. Vieweg \& Sohn, Braunschweig, 1997.

\bibitem{king-schofield}
A.~King and A.~Schofield.
\newblock Rationality of moduli of vector bundles on curves.
\newblock {\em Indag. Math. (N.S.)}, 10(4):519--535, 1999.

\bibitem{laumon}
G.~Laumon and L.~Moret-Bailly.
\newblock {\em Champs alg\'ebriques}.
\newblock Springer-Verlag, Berlin, 2000.

\bibitem{git}
D.~Mumford, J.~Fogarty, and F.~Kirwan.
\newblock {\em Geometric invariant theory}.
\newblock Springer-Verlag, Berlin, 1994.

\bibitem{newstead}
P.~E. Newstead.
\newblock Rationality of moduli spaces of stable bundles.
\newblock {\em Math. Ann.}, 215:251--268, 1975.

\bibitem{newstead2}
P.~E. Newstead.
\newblock Correction to: ``{R}ationality of moduli spaces of stable bundles''.
\newblock {\em Math. Ann.}, 249(3):281--282, 1980.

\bibitem{russo-teixidor}
B.~Russo and M.~Teixidor~i Bigas.
\newblock On a conjecture of {L}ange.
\newblock {\em J. Algebraic Geom.}, 8(3):483--496, 1999.

\bibitem{tyurin}
A.~N. Tyurin.
\newblock Classification of vector bundles over an algebraic curve of arbitrary
  genus.
\newblock {\em Izv. Akad. Nauk SSSR Ser. Mat.}, 29:657--688, 1965.
\newblock English translation: \emph{Amer. Math. Soc. Transl.} 63, 1967.

\bibitem{tyurin2}
A.~N. Tyurin.
\newblock Classification of $n$-dimensional vector bundles over an algebraic
  curve of arbitrary genus.
\newblock {\em Izv. Akad. Nauk SSSR Ser. Mat.}, 30:1353--1366, 1966.
\newblock English translation: \emph{Amer. Math. Soc. Transl.} 73, 1968.

\bibitem{vistoli}
A.~Vistoli.
\newblock Intersection theory on algebraic stacks and on their moduli spaces.
\newblock {\em Invent. Math.}, 97(3):613--670, 1989.

\end{thebibliography}

\end{document}